\documentclass[11pt]{amsart}

\usepackage{amscd}
\usepackage{eucal}
\usepackage{enumerate}
\usepackage{amsthm}
\usepackage{amssymb}
\usepackage{amsmath}
\usepackage{graphics}
\usepackage{epsf}
\usepackage[all]{xy}
\usepackage{psfrag}

\newtheorem{Theorem}{Theorem}[section]
\newtheorem*{theorem*}{Theorem}
\newtheorem{lemma}{Lemma}

\newtheorem*{proposition*}{Proposition}
\newtheorem{proposition}{Proposition}[section]

\newtheorem*{corollary*}{Corollary}
\theoremstyle{definition}

\newtheorem*{comments*}{Comments}
\newtheorem{example}{Example}
\newtheorem*{example*}{Example}

\newtheorem*{remark*}{Remark}
\newtheorem*{remarks*}{Remarks}

\numberwithin{equation}{section}

\gdef\myletter{}

\let\savetheequation\theequation
\def\theequation{\savetheequation\myletter}

\newcommand{\bbC}{{\mathbb C}}

\newcommand{\bbR}{{\mathbb R}}
\newcommand{\bbT}{{\mathbb T}}

\newcommand{\bbZ}{{\mathbb Z}}

\def    \to     {{\longrightarrow}}

\begin{document}

\title{The Fundamental Group of $S^1$-manifolds}
\author{L. Godinho \and M. E. Sousa-Dias}
\thanks{Both authors were partially supported by FCT through program POCTI/FEDER; L. Godinho was partially supported by FCT through grant POCTI/MAT/57888/2004 and by Funda\c{c}\~{a}o Calouste Gulbenkian}
\address{Departamento de Matem\'{a}tica, Instituto Superior T\'{e}cnico, Av. Rovisco Pais, 1049-001 Lisbon, Portugal.}
\email{lgodin@math.ist.utl.pt}
\email{edias@math.ist.utl.pt}
\maketitle
\begin{abstract} We address the problem of computing the fundamental group of a symplectic $S^1$-manifold for non-Hamiltonian actions on compact manifolds, and for Hamiltonian actions on non-compact manifolds with a proper moment map. We generalize known results for compact manifolds equipped with a Hamiltonian $S^1$-action. Several examples are presented to illustrate our  main results.
\end{abstract}
\section{Introduction}
\label{sec:1}
In this paper we address the problem of computing the fundamental group of a symplectic $S^1$-manifold.  For a compact manifold equipped with a Hamiltonian circle action, a result in \cite{L} states that this group is equal to the fundamental group of any of its reduced spaces (as topological spaces) and to the fundamental group of its minimum and maximum level sets. We will consider here  non-Hamiltonian actions on compact manifolds (Theorem~\ref{thm:1})  and Hamiltonian actions on non-compact manifolds with a proper moment map (Theorem~\ref{thm:2}). 

When the action is {\bf non-Hamiltonian}, one can consider a \emph{generalized moment map}\footnote{This generalized moment map is a special case of a Lie group-valued moment map (see \cite{OR} and the references therein).} introduced by McDuff in \cite{MD1} as follows: first, the symplectic form is deformed to a rational invariant symplectic form making the non-zero class $[\iota(\xi_M)\omega]$ rational, where $\xi_M$ denotes the vector field generating the action; then, for a multiple of this symplectic form, there is a map $\phi:M \to S^1$ such that $\iota(\xi_M)\omega = \phi^*(d\theta)$, called  \emph{generalized moment map} (or \emph{circle valued moment map}). This map has many of the properties of an ordinary moment map and can even be used to reduce $M$. In particular, choosing an invariant pair of a Riemannian metric $g$ and a compatible almost complex structure $J$ on $M$ and identifying $S^1$ with $\bbR/\bbZ$ in the usual way, we may define the gradient of $\phi$ with respect to $g$ and see that it is equal to $J\xi_M$.  Its flow has all the nice properties of the gradient flow of an ordinary moment map. In particular, its critical set is a disjoint union of symplectic submanifolds of $M$ (each of codimension at least $4$ since $\phi$ has no local maxima or minima). 

Using the gradient flow of $\phi$ we prove (Theorem~\ref{thm:1}) that, if $M$ is a connected compact symplectic manifold equipped with a non-Hamiltonian circle action and  $P$ is a connected component of an  arbitrary level set of the generalized moment map $\phi$ then, as fundamental groups of topological spaces, $\pi_1(M)$ is a semidirect product
$
\pi_1(M)=\pi_1(P)\rtimes \bbZ,
$ 
when the action has no critical points, or 
$
\pi_1(M)=\pi_1(M_\mathrm{red}) \rtimes \bbZ,
$
where $M_\mathrm{red}:=P/S^1$ is a connected component of the symplectic quotient $\phi^{-1}(a_0)/S^1$ where $a_0=\phi(P)$. 

Note that the proof for the Hamiltonian case presented in \cite{L} relies heavily, at each step, on the existence of a minimum and so it cannot be adapted to the non-Hamiltonian case. Nevertheless, since we have a generalized moment map, we can still use  (circle-valued) Morse theory to prove the above result. 

When the action is Hamiltonian but $M$ is {\bf not compact} one can again use Morse theory, provided that the moment map is proper (i.e. the inverse image of a compact set is compact). In this case we obtain that, if $(M,\omega)$ is a connected symplectic manifold (not necessarily compact) with proper moment map $\phi:M\to \bbR$ and $P$ is an arbitrary (compact) level set of $\phi$ then, as fundamental groups of topological spaces, $\pi_1(M)$ is either,
$\pi_1(M)=\pi_1(P),$
when the action has no critical points, or
$\pi_1(M)=\pi_1(M_\mathrm{red}),$
where $M_\mathrm{red}$ is the symplectic quotient at any value of $\phi$. 
Moreover, if $\phi$ has a local minimum (or a maximum), we recover the referred result for the compact case  in \cite{L}:
$\pi_1(M)=\pi_1(M_\mathrm{red})=\pi_1(F_\mathrm{min})\,\,\mbox{(or equal to $\pi_1(F_\mathrm{max})$)}$,
where $F_\mathrm{min}$ and $F_\mathrm{max}$ are the sets of minimal and maximal points respectively.

Although properness of the moment map is a strong condition which is not verified in many problems in classical mechanics with a global $S^1$-action, our results may still be relevant when, for instance, we can perform a preliminary reduction or symplectic cutting (cf. \cite{Le}) making the induced $S^1$-moment map proper. Let us remark, however,  that the requirement of a proper moment map is essential to our results as can be seen in Examples~\ref{ex:6} and \ref{ex:7}. Indeed, even the statement in Proposition~\ref{prop:3} that all reduced spaces have the same fundamental group may fail to hold when the moment map is not proper. 

Finally, in Section~\ref{sec:5}, we present several other examples illustrating our results.

\section{Non-Hamiltonian circle actions}
\label{sec:2}
In this section we prove our result for non-Hamiltonian actions on compact manifolds:
\begin{Theorem}
\label{thm:1}
Let $M$ be a connected, compact symplectic manifold equipped with a non-Hamiltonian circle action and $\phi:M \to S^1$ the corresponding  generalized moment map. Let $P$ be a connected component of an  arbitrary level set of $\phi$. Then, as fundamental groups of topological spaces, $\pi_1(M)$ is  a semidirect product 
$$
\pi_1(M)=\pi_1(P)\rtimes \bbZ,
$$ 
when the action has no critical points, or
$$
\pi_1(M)=\pi_1(M_\mathrm{red}) \rtimes \bbZ,
$$
where $M_\mathrm{red}:=P/S^1$ is a connected component of the (arbitrary) symplectic quotient $\phi^{-1}(a_0)/S^1$, for $a_0=\phi(P)$.
\end{Theorem}
 Throughout, we  shall choose an $S^1$-invariant compatible pair, $(J,g)$, of an almost complex structure and a Riemannian metric and we  identify $S^1$ with $\bbR/\bbZ$ in the usual way to define the gradient of $\phi$ with respect to $g$. This gradient  is equal to $J\xi_M$, where $\xi_M$ is the vector field generating the action.
 
In order to prove  Theorem~\ref{thm:1}  we will need a series of preliminary results. The first one is proved in \cite{O} but we include a sketch of its proof for the sake of completion.  
\begin{lemma}\cite{O}
\label{le:1}
Let $M$ be a symplectic compact connected manifold equipped with a non-Hamiltonian circle action and  $\phi:M\to S^1$  the corresponding generalized moment map. Then, given any point $y_0\in M$, there exists a homologically non-trivial loop $\gamma : S^1 \to M$ passing through $y_0$.  
\end{lemma}
\begin{proof} 
Since the generalized moment map is locally a function, we can define its Hessian at critical points, their indices and  the gradient flow of $\phi$. Moreover, since the action is non-Hamiltonian, the critical points cannot have index $0$ nor $2n$, where $2n$ is the dimension of $M$. Let us consider the quotient space $X=M/\!\sim$,  where $x\sim y$ iff $x$ and $y$ are in the same connected component of a level set of $\phi$.  As the indices of  the  critical points of $\phi$ are even, $X$ has no branch point. Moreover, $X$ has no boundary and is homeomorphic to a circle (cf. \cite{O} for details). Therefore, we can deform the trajectory of the gradient flow of $\phi$ passing through $y_0$ to a homologically non-trivial loop $\gamma: S^1 \to M$ through $y_0$.
\end{proof}
\begin{lemma}
\label{le:2}
Let $M$ be a compact symplectic manifold equipped with a non-Hamiltonian circle action and   $\phi:M\to S^1$  its generalized moment map. Then, for any regular value $a_0$ and a point $y_0\in \phi^{-1}(a_0)$, the  inclusion  $j:\phi^{-1}(a_0) \to M$ induces an exact sequence of fundamental groups
$$
\pi_1(\phi^{-1}(a_0),y_0) \stackrel{j_*}{\to} \pi_1(M,y_0) \stackrel{\phi_*}{\to} \pi_1(S^1,a_0).
$$
\end{lemma}
\begin{proof}
Clearly $\mathrm{im}(j_*) \subset \ker (\phi_*)$, so we just need to show that  $\ker (\phi_*)\subset \mathrm{im}(j_*)$. Let $[\gamma]\in \ker (\phi_*)$. Then, identifying $S^1$ with $\bbR/\bbZ$, we may assume without loss of generality that there are regular values of $\phi$, $a$ and $b$ with $0\leq a \leq a_0 \leq b <1$, for which $\gamma$ is homotopic to a loop contained in $M^{[a,b]}:=\{x\in M: \, a\leq \phi(x) \leq b\}$. Let  $M^{[a,b]}_{y_0}$ be the connected component of  $M^{[a,b]}$ containing $y_0$. 

If there are no critical points in $M^{[a,b]}_{y_0}$ then $\pi_1(M^{[a,b]},y_0)=\pi_1(\phi^{-1}(a_0),y_0)$  (cf. \cite{Mi}) and so $[\gamma]\in \mathrm{im}(j_*)$. 

If there is just  one critical value $c$ in $(a,b)$ let us consider  $F$, a connected component of the corresponding critical set inside $M^{[a,b]}_{y_0}$ (if there is more than one component we argue similarly for each one). The normal bundle of $F$ has a complex structure induced by the almost complex structure $J$ and splits as a sum $\nu^-\oplus \nu^+$, where $\nu^-$ is tangent to the incoming flow lines of $J\xi_M$ (that is, tangent to the stable manifold). Let $D_F^-$ be the negative disk bundle of $\nu^-$ and $S(D_F^-)$ its sphere bundle.  By Morse theory (see \cite{Mi}) we have
$$
M^{[a,b]}_{y_0}=M^{[a,\tilde{a}]}_{y_0} \cup_{S(D_F^-)} D_F^-,
$$
for any regular value $\tilde{a}$ in $(a,c)$. Hence, by the Van-Kampen theorem, $\pi_1(M^{[a,b]}_{y_0})$ is the free product with amalgamation\footnote{The term amalgamation in $G_1 *_A G_2$ is usually used for the quotient group of the free product of $G_1$ by $G_2$ obtained by identifying the two subgroups that correspond to $A$ under two monomorphisms $A\to G_i$ (see for example \cite{CGKZ}). Here we slightly  abuse this notation since we do not require these maps to be one-to-one.}
\begin{equation}
\label{eq:1}
\begin{array}{ll}
\pi_1(M^{[a,b]}_{y_0})&=\pi_1(M^{[a,\tilde{a}]}_{y_0}) *_{\pi_1(S(D_F^-))} \pi_1(D_F^-) \\
&\\
&= \pi_1(\phi^{-1}(\tilde{a}), \tilde{y}) *_{\pi_1(S(D_F^-))} \pi_1(F),
\end{array}
\end{equation}
where $\tilde{y}$ is a point in the appropriate component of $\phi^{-1}(\tilde{a})$.

If $\mathrm{index}(F) > 2$, then $\pi_1(S(D_F^-))$ is isomorphic to $\pi_1(F)$ and so, since we also have $\pi_1(D_F^-)=\pi_1(F)$, we get $\pi_1(M^{[a,b]}_{y_0}) = \pi_1(\phi^{-1}(\tilde{a}),\tilde{y})$.

If $\mathrm{index}(F)=2$, we consider the principal circle bundle
$$
S^1 \stackrel{\tilde{i}}{\hookrightarrow} S(D_F^-) \stackrel{\tilde{p}_S}{\to} F
$$
and its homotopy exact sequence
$$
\cdots \to \pi_1(S^1) \stackrel{\tilde{i}_*}{\to} \pi_1(S(D_F^-))\stackrel{(\tilde{p}_S)_*}{\to} \pi_1(F) \to \{1\}.
$$
Note that $S(D_F^-)$ can be identified with the restriction of the circle bundle $\phi^{-1}(\tilde{a})\to M_{\tilde{a}}$ to $F$, where $M_{\tilde{a}}$ is the reduced space $\phi^{-1}(\tilde{a})/S^1$ (there is an embedding of $F$ in $M_{\tilde{a}}$ as it is shown in \cite{L}), and so we also have an inclusion $S(D_F^-)\stackrel{\kappa}{\hookrightarrow}\phi^{-1}(\tilde{a})$.
In the amalgamation  (\ref{eq:1}), the elements of $(\tilde{p}_S)_*(\pi_1(S(D_F^-)))=\pi_1(F)$ (the map $(\tilde{p}_S)_*$ is surjective) are identified with the corresponding elements in $\kappa_*(\pi_1(S(D_F^-)))\subset \pi_1(\phi^{-1}(\tilde{a}),\tilde{y})$, implying that $\pi_1(M^{[a,b]}_{y_0})$ can be identified with the quotient $\pi_1(\phi^{-1}(\tilde{a}),\tilde{y})/N$, where $N$ is the normal subgroup generated by all the elements of $\kappa_*(\ker{((\tilde{p}_S)_*)})$. 

Repeating this argument using $-\phi$ instead of $\phi$ (corresponding to reversing  the direction of  the circle action), we can substitute $\tilde{a}$ by any value $\hat{a}\in (c,b]$ in the above argument. However, the relevant critical points will no longer be the index-$2$ critical points but the ones with index equal to $2n-2$, where $2n$ is the dimension of $M$. 

If $[a,b]$ has more than one critical value, let $n_1$ be the number of critical values in $(a_0,b]$ for which there is an  index-$2$ component of the corresponding critical set intersecting $M^{[a,b]}_{y_0}$. Similarly, let $n_2$ be the number of critical values in $[a,a_0)$ for which there is an  index-$(2n-2)$ component of the corresponding critical set intersecting $M^{[a,b]}_{y_0}$. By induction on $n_1$ and $n_2$ and by using the Van-Kampen Theorem as in (\ref{eq:1}) each time we cross one of those critical levels (using $\phi$ or $-\phi$ accordingly) we see that $\pi_1(M^{[a,b]}_{y_0},y_0)$ can be obtained from $\pi_1(\phi^{-1}(a_0),y_0)$ by taking a sequence of $n_1+n_2$ quotients as it is explained above, and the result follows.
\end{proof}
Note that the level sets of the generalized moment map $\phi$ may not be connected leading to non-connected reduced spaces (cf. Example~\ref{ex:1}). Nevertheless, we will show that all their connected components have the same fundamental group. For that, we first consider the equivalence relation $\sim$ defined in the proof of  Lemma~\ref{le:1} and take $M/\!\!\sim \, \cong S^1$. The map $\phi$ descends to the quotient $M/\!\!\sim$ (since $y\sim x$ implies $\phi(x)=\phi(y)$) giving us a finite covering of $S^1$, $\tilde{\phi}:M/\!\!\sim\, \cong S^1 \to S^1$ (i.e. $\tilde{\phi}(z)=z^k$ for $z\in S^1$ and some $k\in \bbZ$). Hence, we have the following decomposition of $\phi$
\begin{equation}\label{eq:1.1}
\xymatrix{
M\ar[r]^(0.4){\phi^c} \ar@/_1.5pc/[rr]^\phi &M/\!\!\sim\, \cong S^1\ar[r]^(0.6){\tilde{\phi}} & S^1,
}
\end{equation}
where the map $\phi^c:M\to S^1$ is surjective and has connected level sets. Moreover, considering the gradient flow of $\phi^c$ with respect to the metric $g$, it is easy to check that it has the same critical set  as $\phi$ as well as all the nice properties of its gradient flow. In particular, the indices of the critical submanifolds are all even. Moreover, the $k$ connected components of the reduced space $M_a:=\phi^{-1}(a)/S^1$ are the reduced spaces $M_{a_j}^c:=(\phi^c)^{-1}(a_j)/S^1$ of $\phi^c$, where the $a_j$'s ($j=1,\ldots,k$) are such that $\tilde{\phi}(a_j)=a$, that is,  $M_a$ is the disjoint union
$$
M_a= \bigsqcup_{j=1}^k (\phi^c)^{-1}(a_j)/S^1.
$$
\begin{proposition}
\label{prop:3}
Let $M$ be a manifold satisfying the hypotheses of Lemma~\ref{le:1}. Then, the fundamental group of all connected components of all reduced spaces  $M_a:=\phi^{-1}(a)/S^1$ is always the same, even for critical values of the generalized moment map.
\end{proposition}
\begin{proof} 
Let us consider the map $\phi^c: M\to S^1$ defined above.
If the action has no fixed points then all ``reduced spaces'' $(\phi^c)^{-1}(a)/S^1$ are diffeomorphic and we are done. 
If that is not the case, let us again identify $S^1$ with $\bbR/\bbZ$  and assume that $0$ is a regular value of $\phi^c$ (if not, we just break up the circle at another point). Let $c_1$ be the smallest critical value of $\phi^c$ in $[0,1]$ and consider a connected component $F$ of the corresponding critical set. Let $D_F^-$ be the negative disk bundle of $\nu^-$ and $S(D_F^-)$ its sphere bundle. Then, by Morse theory, $(\phi^c)^{-1}(c_1)$   has the same homotopy type as $(\phi^c)^{-1}(a)\cup_{S(D_F^-)} D_F^-$ where $a$ is any regular value in $[0,c_1)$. This implies that $M_{c_1}^c:=(\phi^c)^{-1}(c_1)/S^1$ has the same  homotopy type as $\left((\phi^c)^{-1}(a)/S^1\right)\cup_{S(D_F^-)/S^1} \left(D_F^-/S^1 \right) = M_a^c \cup_{S(D_F^-)/S^1}\left(D_F^-/S^1\right)$, and so $\pi_1(M_{c_1}^c)$ is the free product with amalgamation
$$
\pi_1(M_{c_1}^c)=\pi_1(M_a^c) *_{\pi_1(S(D_F^-)/S^1)} \, \pi_1(D_F^-/S^1).
$$ 
However, the local normal form for $\phi$ (and consequently for $\phi^c$) on a neighborhood of $F$ is the same as the local normal form of a neighborhood of a critical set of an ordinary moment map, implying that $S(D_F^-)/S^1$ is a weighted projectivized bundle over $F$ and so, since we also have that $D_F^-$ is homotopy equivalent to $F$, we conclude that 
$$\pi_1(S(D_F^-)/S^1)=\pi_1(D_F^-/S^1)$$
and so $\pi_1(M_{c_1}^c)=\pi_1(M_a^c)$. Using $-\phi$ instead of $\phi$ (corresponding to reversing the direction of the circle action) we obtain that $\pi_1(M_{c_1}^c)=\pi_1(M_b^c)$ for $b\in (c_1,c_2)$ where $c_2$ is a critical value of $\phi^c$ (if it exists) and the interval $(c_1,c_2)$  contains only regular values. Repeating this for every critical value of $\phi^c$ we conclude that all connected components of all reduced spaces (even critical ones) have the same fundamental group. 
\end{proof}
Note that an alternative proof of this Proposition could follow from the fact that, for circle actions, when passing a critical value of the moment map the reduced spaces change by a weighted blow-down followed by a weighted blow-up (cf. \cite{G} and \cite{BP}) and so, by \cite{MD2}, their fundamental group does not change.
\begin{lemma}
\label{le:4}
Let $M$ be a manifold satisfying the hypotheses of Lemma~\ref{le:1} equipped with a circle action with a non-empty fixed point set. Let $a_0 \in S^1$ be a  regular value of the generalized moment map $\phi$, and consider the principal circle bundle 
$$S^1 \stackrel{i}{\hookrightarrow} \phi^{-1}(a_0) \stackrel{p}{\to} M_{a_0},$$
where $M_{a_0}:=\phi^{-1}(a_0)/S^1$ is the reduced space at $a_0$. Then, for $y_0\in \phi^{-1}(a_0)$,  the kernel of the map $p_*:\pi_1(\phi^{-1}(a_0),y_0)\to \pi_1(M_{a_0}, p(y_0))$ is equal to the kernel of the map $j_*:\pi_1(\phi^{-1}(a_0),y_0)\to \pi_1(M, y_0)$ defined in Lemma~\ref{le:2}.
\end{lemma}
\begin{proof}
Clearly $\ker{(p_*)}\subset \ker{(j_*)}$. Indeed, if $[\gamma]\in \ker{(p_*)}$ then, the gradient flow of $\phi$ gives us a homotopy between $\gamma$ and a constant path contained in some critical level set and so $[\gamma]=1$ in $\pi_1(M)$.

Let us now see that $\ker{(j_*)}\subset \ker{(p_*)}$. Take $[\gamma]\in \ker{(j_*)}$. Then, $\gamma$ is homotopic to a nullhomotopic loop in $M^{[a,b]}_{y_0}$ for some regular values $a,b$ with $0\leq a \leq a_0 \leq b <1$. Indeed, if that were not the case, there would exist a homotopy $H:[0,1]\times [0,1] \to M$ between $\gamma$ and the constant path based at $y_0$ for which $\phi$, restricted to $D:= \mathrm{im}(H)$ would be surjective. Hence, there would exist a loop $\alpha:S^1\to D$ in $D$ such that $(\phi\circ \alpha) (S^1)=S^1$ and then, since $\pi_1(D)=\{1\}$, we would have $[\alpha]=1$ in $\pi_1(D)$ while $\phi_*[\alpha]=[\phi\circ\alpha]\neq 1$, which is impossible.  

If the critical points in $M^{[a,b]}_{y_0}$ have index greater than $2$ and smaller than $2n-2$ (where $2n$ is the dimension of $M$), or if there are no critical points at all in this set, then, as we saw in the proof of  Lemma~\ref{le:2}, $\pi_1(M^{[a,b]}_{y_0},y_0)=\pi_1(\phi^{-1}(a_0),y_0)$, implying that $\gamma$ is nullhomotopic in $\phi^{-1}(a_0)$ and so $p_*([\gamma])=1$. 

If there are components of the critical set  inside $M^{[a,b]}_{y_0}$  with index equal to $2$ or equal to $2n-2$  then, again like in  the proof of Lemma~\ref{le:2}, $\pi_1(M^{[a,b]}_{y_0},y_0)$ can be obtained from $\pi_1(\phi^{-1}(a_0),y_0)$ by taking a sequence of quotients. Indeed, for each index-$2$ component $F$ with $\phi(F)=c_i > a_0$, we consider the maps $\kappa_i:S(D_F^-)\to \phi^{-1}(a_0)$, $p_{S_i}:S(D^-_F) \to F$ and $p:\phi^{-1}(a_0)\to M_{a_0}$ defined as in the proof of Lemma~\ref{le:2}; then we take a sequence of quotients of  $\pi_1(\phi^{-1}(a_0),y_0)$ by $N_i$, the normal subgroups generated by all the elements of $(\kappa_i)_*(\ker{(p_{S_i})_*)}$; we repeat this procedure for each index-$(2n-2)$ component $F$ with $\phi(F) < a_0$, this time using $-\phi$ instead of $\phi$. We conclude that, if $[\gamma]=1$ in $\pi_1(M^{[a,b]}_{y_0},y_0)$,  then  $[\gamma]\in N_i$ for one of the groups $N_i$ considered above. However, $(\kappa_i)_*(\ker{(p_{S_i})_*)}\subset \ker{(p_*)}$ and so, since $\ker{(p_*)}$ is normal, we conclude that $[\gamma]\in N_i \subset \ker{(p_*)}$ and the result follows.
\end{proof}
With these results  we can now prove Theorem~\ref{thm:1}. 
\begin{proof}({\it of Theorem~\ref{thm:1}})
First, let us assume that the action has no fixed points that is, the generalized moment map $\phi$ has no critical points. In this case, since $M$ is connected, all the level sets of $\phi$ are equivariantly diffeomorphic, since we can  use the flow of $J\xi_M$ to identify the level sets. Moreover, since we are assuming that there are no fixed points, the map $\phi^c$ defined in (\ref{eq:1.1}) is a fibration with connected fiber  $P:=(\phi^c)^{-1}(a)$ (a fixed level set of $\phi^c$) which is a connected component of the level set $\phi^{-1}(a^k)$ for some $k\in \bbZ$. Hence, the long exact homotopy sequence for $P \to M \to S^1$ gives us that the sequence
\begin{equation}
\{1\} \to \pi_1(P) \stackrel{j_*}{\to} \pi_1(M) \stackrel{(\phi^c)_*}{\to} \pi_1(S^1)\to \{1\} 
\end{equation}
is exact, implying that $j_*$ is injective. Moreover, the homologically non-trivial loop $\gamma:S^1\to M$ given by Lemma~\ref{le:1} is a section of the above fibration, and so $\phi^c_*\circ \gamma_* =\mathrm{id}$.    Hence,  $\pi_1(M)$ is a semidirect product  $\pi_1(M) = G_1 \rtimes G_2$ where $G_1$ is the kernel of $\phi^c_*$ and $G_2:=\mathrm{im}(\gamma_*)$.  Moreover, we have $G_1=\mathrm{im}(j_*) \cong \pi_1((\phi^c)^{-1}(a))/\ker{(j_*)}=\pi_1((\phi^c)^{-1}(a))$ (since  $j_*$ is injective), and so, as $(\phi^c)_*$ maps $G_2$ isomorphically onto $\pi_1(S^1)=\bbZ$, the result follows.

If the action has fixed points they cannot be local maxima nor minima. Taking a fixed point $F$, we consider the homologically non-trivial loop $\gamma: S^1 \to M$ through $F$ given by Lemma~\ref{le:1}. 
By Lemmas~\ref{le:2} and \ref{le:4}, we have the following exact sequence, where $a$ is a regular value of $\phi^c$:
\begin{equation}\label{eq:2}
\xymatrix{
\pi_1(S^1)\ar[r]^(0.4){i_*} &\pi_1\left(\left(\phi^c\right)^{-1}(a)\right)\ar[r]^(0.6){j_*} &\pi_1(M)\ar[r]^{(\phi)_*^c} &\pi_1(S^1)\ar@/^1pc/[l]^{\gamma_*}
}
\end{equation}
with $(\phi)^c_*\circ \gamma_*=\mathrm{id}$.
Hence, taking  $G_1:=\ker{((\phi^c)_*)}$ and $G_2:=\mathrm{im}(\gamma_*)$ we have that $\pi_1(M)$ is a  semidirect product of $G_1$ and $G_2$. Moreover, considering the map  $p: (\phi^c)^{-1}(a) \to (\phi^c)^{-1}(a)/S^1=: M_a^c$, we have, by Lemma~\ref{le:4}, that  $\ker{(j_*)}=\ker{(p_*)}$, and so 
$$G_1=\mathrm{im}(j_*)= \pi_1((\phi^c)^{-1}(a))/\ker{(p_*)}= \pi_1(M^c_a),$$
where the  ``reduced space"  $M^c_a$   is a connected component of the reduced space $\phi^{-1}(a^k)/S^1$ of $\phi$. On the other hand, $\phi^c_*$ maps $G_2$ isomorphically onto $\pi_1(S^1)=\bbZ$. Hence, $\pi_1(M)$ contains two subgroups  $G_1$ and $G_2$ such that $G_1$ is normal and isomorphic to $\pi_1(M^c_a)$ and $G_2$ is isomorphic to $\bbZ$. Moreover, each element of $\pi_1(M)$ is uniquely represented as the product of an element of $G_1$ by an element of $G_2$. Indeed, $\pi_1(M)$ is a semidirect product of $\pi_1(M^c_a)$ by $\bbZ$ and then, by Proposition~\ref{prop:3}, the result follows.
\end{proof}
\begin{remark*}
To be able to completely determine the semidirect product above one must know how the elements of $\bbZ$  ``act'' by conjugation on the fundamental group of a connected component $P$ of a level set of the moment map (when the action has no fixed points) or on the fundamental group of a connected component of a reduced space, $M_\mathrm{red}=P/S^1$. Indeed, one needs to establish the homomorphism $\Psi:\bbZ \to Aut(\pi_1(P))$ or $\Psi:\bbZ \to Aut(\pi_1(M_\mathrm{red}))$ given by $\Psi(j)(g)=j g j^{-1}$, where $j\in \bbZ$ (since $\bbZ$ is cyclic it suffices to know the image of the generator). This will of course  depend on the manifold $M$. Nevertheless, this   ``action'' of $\bbZ$ is independent of the choice of the level set $P$.
\end{remark*}
\section{Non-Compact Manifolds}
\label{sec:3}
We consider in this section Hamiltonian circle actions with a proper moment map on non-compact manifolds $M$. The proof that the fundamental group of $M$ is equal to the fundamental group of its reduced spaces does not follow from the proof for the compact case presented in \cite{L} since, in this case, we do not necessarily have a maximum or a minimum.   
\begin{Theorem}
\label{thm:2}
Let $S^1$ act on a connected symplectic manifold $(M,\omega)$ (not necessarily compact) with proper moment map $\phi:M\to \bbR$, and let $P$ be an arbitrary (compact) level set of $\phi$. Then, as fundamental groups of topological spaces, the fundamental group of $M$ is either,
$$\pi_1(M)=\pi_1(P),$$
when the action has no critical points, or
$$\pi_1(M)=\pi_1(M_\mathrm{red}),$$
where $M_\mathrm{red}$ is the symplectic quotient at any value of $\phi$. 

Moreover, if $\phi$ has a local minimum (or a maximum), then 
$$\pi_1(M)=\pi_1(M_\mathrm{red})=\pi_1(F_\mathrm{min})\,\,\mbox{(or $=\pi_1(F_\mathrm{max})$)},$$
where $F_\mathrm{min}$ and $F_\mathrm{max}$ are the sets of minimal and maximal points respectively.
\end{Theorem} 
\begin{proof} Just as in the compact case, the image of the moment map is an interval $I\subset \bbR$ but now not necessarily compact. However,  the level sets of $\phi$ are still connected \cite{LMTW}.
If $\phi$ has no critical points, then, as in classical Morse theory, $M$ is diffeomorphic to $\phi^{-1}(a)\times I$ for any value $a$ of $\phi$, and so $\pi_1(M)=\pi_1(\phi^{-1}(a))$. 

If the action has fixed points, then, considering  a regular value of $\phi$, $a_0$, we can adapt the proof of Lemma~\ref{le:2} to show that the sequence
$$
\pi_1(\phi^{-1}(a_0))\stackrel{j_*}{\hookrightarrow}\pi_1(M) \stackrel{\phi_*}{\to} \{1\}
$$
is exact (i.e. $j_*$ is surjective). Indeed, if $[\gamma]\in \pi_1(M)$, then $\gamma$ is homotopic to some loop contained in a compact set ($\phi$ is proper)  $M^{[a,b]}:=\{x\in M:\, a\leq \phi(x) \leq b\}$ for some values $a,b\in \bbR$ with $a\leq b$. 

The proof of Proposition~\ref{prop:3} can also be adapted to show that all reduced spaces (even critical ones) have the same fundamental group. Similarly, we can use  the proof of Lemma~\ref{le:4} to show that the kernel of the map $p_*:\pi_1(\phi^{-1}(a_0))\to \pi_1(M_{a_0})$ is equal to the kernel of the map $j_*:\pi_1(\phi^{-1}(a_0))\to \pi_1(M)$, induced respectively by the quotient and the inclusion maps (here $M_{a_0}$ denotes the reduced space $\phi^{-1}(a_0)/S^1$). Consequently, since $j_*$ and $p_*$ are surjective,
$$
\pi_1(M_{a_0})=\pi_1(\phi^{-1}(a_0))/\ker{(p_*)}=\pi_1(\phi^{-1}(a_0))/\ker{(j_*)}=\pi_1(M).
$$    
Hence, to finish our proof, we just need to show that, when $\phi$ has either a local minimum at $F_\mathrm{min}$, or a local maximum at $F_\mathrm{max}$, we have $\pi_1(M_{a_0})=\pi_1(F_\mathrm{min})$ or $\pi_1(M_{a_0})=\pi_1(F_\mathrm{max})$. Let us consider the case where $\phi$ has a minimum (the other case is similar). Here we can use the following argument used in \cite{L} for the compact case: let $m$ be the minimum value of $\phi$ and consider an interval $(m,b)$ formed by regular values of $\phi$. For $a\in (m,b)$ we have, by the equivariant symplectic embedding theorem, that $\phi^{-1}(a)$ is a sphere bundle over $F_\mathrm{min}$. Let $S^{2l+1}$ be its fiber, where $\mathrm{dim} (F_\mathrm{min})=2(n-l-1)$. Then, the reduced space $M_a$ is diffeomorphic to an orbibundle over $F_\mathrm{min}$  with fiber a weighted projective space $\bbC P^{l}:=S^{2l+1}/S^1$, and we have the exact sequence
$$
\pi_1(\bbC P^l) \to \pi_1(M_a) \rightarrow \pi_1(F_\mathrm{min})\to \{1\}.
$$     
Since $\bbC P^l$ is simply connected, we have $\pi_1(F_\mathrm{min})=\pi_1(M_a)$ and the result follows.
\end{proof}
\newpage
\section{Examples}
\label{sec:5}
\subsection{Non-Hamiltonian actions on compact manifolds}\hspace{1mm}

\begin{example} 
\label{ex:1}
Let us begin with a very simple example of a non-Hamiltonian circle action with an empty fixed point set. Let $M$ be the $2$-torus $\bbT^2=S^1 \times S^1$ ($\pi_1(M)=\bbZ^2$) with symplectic form $\sigma=d\theta_1 \wedge d\theta_2$, and consider the  $S^1$-action given by $e^{i\beta} \cdot (e^{i\theta_1},e^{i\theta_2})=(e^{i(2\beta+\theta_1)}, e^{i\theta_2})$. The generalized moment map $\phi:M\to S^1$ is just $\phi(e^{i\theta_1},e^{i\theta_2})=e^{2i\theta_2}$. All the level sets of $\phi$ are equal to two disjoint copies of  $S^1$. We can decompose $\phi=\tilde{\phi} \circ \phi^c$ in the following way:
$$
\begin{array}{ccccc}
M                              &  \stackrel{\phi^c}{\to} & S^1           & \stackrel{\tilde{\phi}}{\to} & S^1 \\
(e^{i\theta_1},e^{i\theta_2})   &  \mapsto                & e^{i\theta_2} & \mapsto                      & e^{2i\theta_2},
\end{array}
$$ 
where the level sets of $\phi^c$ are the connected components of the level sets of $\phi$ (they are all equal to $S^1$), and we get the result in Theorem~\ref{thm:1}, that is, $\pi_1(\bbT^2)$ is a semidirect product of $\pi_1((\phi^c)^{-1}(a),y_0)$ by $\bbZ$. Indeed, for a point  $y_0\in M$ and $a=\phi^c(y_0)$, the two subgroups of $\pi_1(M,y_0)$ isomorphic to $\pi_1((\phi^c)^{-1}(a),y_0)=\bbZ$ and to $\pi_1(S^1)=\bbZ$ are both normal, implying that  their semidirect product is just the regular direct product of the two groups.
\end{example}

\begin{example} 
\label{ex:2}
Let us consider the example of a $6$-manifold $M$ with a free symplectic circle action with {\bf contractible orbits} constructed in \cite{K}. Here, we take $Y$, the smooth oriented simply-connected $4$-manifold underlying a $K3$ surface (see for example \cite{BHPV}), and consider the mapping torus\footnote{The mapping torus of a map $h:Y\to Y$ is the identification space $T(h)=Y\times [0,1]/\{(x,0)=(h(x),1))\mid \, x\in Y\}$ (cf. \cite{R}).} $X$ of an orientation-preserving diffeomorphism $\Phi:Y \to Y$ obtained as follows: first, knowing that the intersection form of $Y$ is $Q=3H\oplus 2E_8$ \footnote{Here $H$ is the hyperbolic plane and $E_8$ is the unimodular even positive definite form of rank $8$ (see \cite{BHPV} for details).}, we consider an automorphism $f$ of $H\oplus H$, such that $f(x)=x+c$ and $f(c)=c$, where $x$ and $c$ are two non-zero primitive classes in $H^2(Y,\bbZ)$; then, we extend $f$ to all of $Q$, preserving the orientation of a maximal positive-definite subspace and find, by a result of Matumoto \cite{M},  an orientation-preserving diffeomorphism $\Phi:Y\to Y$ with $\Phi^*=f$. Finally, from $X$, we obtain $M$ as the total space of the circle bundle $\pi:M\to X$ with Euler class $c$ (since $\Phi^*c=c$, we can choose a lift to a cohomology class on the mapping torus which we also denote by $c$). Since this bundle has contractible fibers, its homotopy long exact sequence gives us that $\pi_1(M)=\pi_1(X)$. Moreover, since $X$ is a mapping torus and $Y$ is simply connected (implying that $\Phi_*:\pi_1(Y)\to \pi_1(Y)$ is trivially an isomorphism), we have that $\pi_1(X)$, a semidirect product $\pi_1(Y) \rtimes \bbZ$ (see \cite{R}), is equal to $\bbZ$. 

Let us now obtain the same result using Theorem~\ref{thm:1}. Let $\omega$ be the $S^1$-invariant symplectic form in $M$ (we omit its construction for simplicity but the details can be found in \cite{K}) and let $\xi_M$ be the vector field generating the action. Since the closed $1$-form $\iota(\xi_M)\omega$ vanishes on the tangents to the $S^1$-action, it can be written as $\pi^*\alpha$ where $\alpha$ is a closed $1$-form on the quotient $X$. Moreover, it is shown in \cite{FGM} that there is a map $\nu:X\to S^1$ for which  $\alpha=\nu^* (d\theta)$. Hence $\iota(\xi_M)\omega= \pi^* \nu^* (d\theta)$ and so $\phi:=\nu \circ \pi$ is the generalized moment map for this action. However, $X$ is a mapping torus implying that there is a natural map $\nu^c$ from $X$ to $S^1$ with (connected) level sets equal to $Y$ and, as is shown in \cite{FGM}, $\nu =\tilde{\nu} \circ \nu^c$ where the map $\tilde{\nu}:S^1 \to S^1$ is a  finite covering of the circle. Hence, the connected components of the level sets of $\nu$ are equal to $Y$ and so, denoting by $P$ a connected component of an arbitrary level set of $\phi$, we get $\pi_1(P)=\pi_1(Y)=\{1\}$ since both the orbits of the circle action and $Y$ are contractible. Therefore, Theorem~\ref{thm:1} also gives $\pi_1(M,y_0)=\pi_1(P,y_0)\rtimes \bbZ=\bbZ$, where $y_0\in P$.   
\end{example}

\begin{example} \label{ex:3}
The only known example of a manifold $M$ equipped with a   non-Hamiltonian circle action with fixed points  was constructed by  McDuff in \cite{MD1}. Theorem~\ref{thm:1} will allow us to compute its fundamental group. This $6$-dimensional manifold, $M$, is obtained  by first considering a special manifold with boundary, $X$, equipped with a Hamiltonian circle action with moment map $\nu:X \to [0,7]$, having two boundary components (lying over the endpoints $0$ and $7$), and then gluing them together. 

This manifold $X$, which has $4$ critical levels at  $s=1$, $2$, $5$ and $6$ with zero sets of codimension $4$, is constructed as follows (for simplicity we will not keep track of symplectic forms): considering coordinates $\theta_1,\ldots, \theta_4$ on $T^4$ and letting $\sigma_{ij}$ be the form $d\theta_i \wedge d\theta_j$, we construct five regular pieces  $\nu^{-1}(I)$, where
\begin{itemize}
\item[$\bullet$] $\nu^{-1}(I) = T^4\times S^1 \times I$, for $ I=(0,1)$ and $I=(6,7)$,  
\item[$\bullet$]  $\nu^{-1}(I) =  P_I \times I$,  with  $P_I$ a principal circle bundle over $T^4$ of Chern class $c_I$, where
\begin{enumerate}
\item[i)]  $c_I:=-[\sigma_{42}]$ for $I=(1,2)$;
\item[ii)] $c_I:=-[\sigma_{31} + \sigma_{42}]$ for $I=(2,5)$;
\item[iii)] $c_I:=-[\sigma_{31}]$ for $I=(5,6)$.
\end{enumerate}
\end{itemize}
Then, we construct four additional pieces $Q_\lambda$ ($\lambda =1,2,5,6$), lying over the intervals $[\lambda-\varepsilon, \lambda+ \varepsilon]$, which are  then glued  to the already defined parts. The singularity as $s$ increases through $1$ is diffeomorphic to the singularity as $s$ decreases through $6$, and similarly for $2$ and $5$, so $X$ will be completely determined with the description of $Q_1$ and $Q_2$. 

The piece $Q_1$ is of the form $T^2 \times Y$, where $Y$ is a $4$-manifold obtained from $S^2 \times S^2$ with symplectic form $2\rho \oplus \rho$ (where $\rho$ is a symplectic form on $S^2$ with total area $1$) and the standard diagonal circle action with moment map $H$ \footnote{\label{fn:1}The moment map $H$ is $2 \mu_1 + \mu_2$ where $\mu_i$ is the moment map for the $i$-th factor with respect to $\rho$; note that $\mu_i(S^2)=[0,1]$.}, in the following way: taking  $V:=H^{-1}([2-\varepsilon, 2+\varepsilon])$, which fibers over $S^2$, cut the inverse image of a disc avoiding the unique critical value of this projection (which is an $S^1$-invariant set diffeomorphic to $D^2\times S^1 \times [-\varepsilon, \varepsilon]$) and glue back a copy of $(T^2-\mathrm{Int}(D^2))\times S^1 \times [-\varepsilon, \varepsilon]$ (cf. Figure~\ref{fig:1}). The piece $Q_2$ is of the form $S \times_{S^1} Y$ where $S$ is the total space of the  circle bundle $S^1 \hookrightarrow S \to T^2$ with Euler characteristic $-1$, so that $Q_2$ fibers over $T^2$ with fiber $Y$.   
\begin{figure}[h!]
\begin{center}
\epsfxsize=.4\textwidth
    \psfrag{H}{$\scriptsize{H(M)}$}
    \psfrag{0}{$\scriptsize{0}$} 
    \psfrag{1}{$\scriptsize{1}$}
    \psfrag{2}{$\scriptsize{2}$}
    \psfrag{3}{$\scriptsize{3}$}
    \leavevmode
    \epsfbox{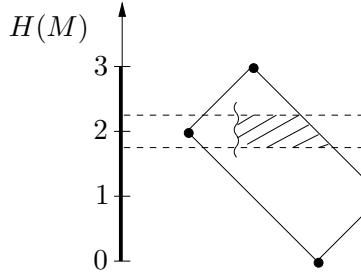}
\end{center}
\caption{Obtaining $Y$ from $S^2\times S^2$ in Example~\ref{ex:3}.}\label{fig:1}
\end{figure}

The manifold $M$ is then obtained from $X$ by gluing $\nu^{-1}(0)$ to $\nu^{-1}(7)$ by the diffeomorphism of $T^4$ that interchanges $\theta_1$ with $\theta_3$ and $\theta_2$ with $\theta_4$.

Any reduced space $M_a$ at a regular value $a$ of the generalized moment map is diffeomorphic to $T^4$, implying that $\pi_1(M_a)=\bbZ^4$, and so, by Proposition~\ref{prop:3} $\pi_1(M_c)=\bbZ^4$, for every reduced space at a critical value $c$. We conclude from Theorem~\ref{thm:1} that $\pi_1(M)=\bbZ^4 \rtimes \bbZ$. The  ``action'' of $\bbZ$ on $\bbZ^4$ is determined by the diffeomorphism of $T^4$ used to glue the boundary components of $X$.
\end{example}

\subsection{Hamiltonian actions on non-compact manifolds}\hspace{1mm}

The first two examples below satisfy the hypotheses of Theorem~\ref{thm:2}. On Example~\ref{ex:4}, the proper moment map has no minima nor maxima while, on Example~\ref{ex:5}, such type of critical points do exist. The last two examples (\ref{ex:6} and \ref{ex:7}) illustrate that the  properness of the   moment map is essential to our results on the fundamental group. In particular, in Example~\ref{ex:6}, there are no critical points and $\pi_1(M)\neq \pi_1(\phi^{-1}(a))$ for some values $a$ of the moment map $\phi$ and, in Example~\ref{ex:7}, there is a critical point (a minimum) and $\pi_1(M_{\mathrm{red}})$ is not always the same for all values of the moment map.
\begin{example} 
\label{ex:4}
We can construct a non-compact symplectic manifold $X$ with a Hamiltonian  $S^1$-action with no minima or maxima from Example~\ref{ex:3} above in the following way: taking the manifold $X$ from McDuff's example we attach two pieces to its boundary of the form $T^4\times S^1 \times I$ where $I=(-\infty, 0]$ and $[7, \infty)$, extending its symplectic form and moment map $\nu$ in the natural way. The resulting  moment map is proper and has no minimum nor maximum. Then, since the fundamental group of its reduced spaces is $\bbZ^4$, so is $\pi_1(X)$.    
\end{example}

\begin{example} 
\label{ex:5}
Consider  $M=S^2\times \bbR^2$ with symplectic form $\omega=\rho_0 \oplus \omega_0$ (where $\rho_0$ and $\omega_0$ are the standard symplectic forms on $S^2$ and on $\bbR^2$), and  the following $S^1$-action: take the $S^1$-action on $S^2$ by rotations about the vertical $z$-axis and the standard $S^1$-action on $\bbR^2$ by rotations around the origin. The moment map on $M$ is just the sum of the  height function $z$ with the  map $\nu (u,v)=(u^2+v^2)/2$. Physically, we have a classical spin and a harmonic oscillator $\nu$. 

This moment  map $\phi=\nu + z$ is proper, has a minimum at $S\times \{0\}$ and a critical point of index $2$ at $N\times \{0\}$, where $S$ and $N$ are respectively the south and north poles of the sphere. This circle action extends to a Hamiltonian $\bbT^2=S^1\times S^1$ torus action, where the action of the second circle  on the sphere  is by clockwise rotations  and on  $\bbR^2$   is the standard  one. The moment map for this extended $\bbT^2$-action is $(\nu +z, \nu -z)$ and its image is pictured in  Figure~\ref{fig:2}. All regular reduced spaces of $\phi$ are homeomorphic to $S^2$ and so, $\pi_1(M)=\pi_1(M_\mathrm{red})=\pi_1(F_\mathrm{min})=\{1\}$.
\begin{figure}[h!]
\begin{center}
\epsfxsize=.4\textwidth
    \psfrag{H}{$\scriptsize{\phi(M)}$}
    \psfrag{-1}{$\scriptsize{-1}$} 
    \psfrag{0}{$\scriptsize{0}$}
    \psfrag{1}{$\scriptsize{1}$}
    \leavevmode
    \epsfbox{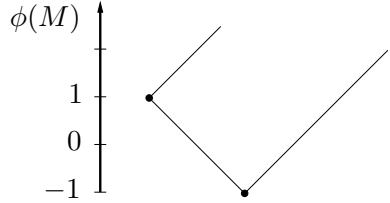}
\end{center}
\caption{Moment map image for the $\bbT^2$-action on $S^2\times \bbR^2$.}\label{fig:2}
\end{figure}
\end{example}

\begin{example}
\label{ex:6}
 Let us now give an example which shows  that  the requirement for properness of the moment map is essential to our results. Let us consider $M=S^2\times \bbC \setminus \{0\}$ with the same symplectic form and the same circle action as in Example~\ref{ex:5} above. The image of the  moment map  $\phi=\nu + z$   is now the interval $(-1,\infty)$; this map {\bf has no critical points} on $M$ and is {\bf no longer proper} (note for instance that the level sets $\phi^{-1}(a)$ for values $a\in (-1,1)$ are not compact). We can easily check that Theorem~\ref{thm:2} is no longer valid. In fact, the fundamental group of the manifold, $\pi_1(M)=\bbZ$, is no longer the fundamental group of the level sets of $\phi$. Indeed, considering, for instance, a value $a\in (-1,1)$, the level sets $\phi^{-1}(a)$ are diffeomorphic to $S^3\setminus\{pt\}$. Note also that the level sets of this moment map are no longer all diffeomorphic (as they would be for a proper moment map with no critical points) since, for $a\in (1,\infty)$, they are diffeomorphic to $S^2\times S^1$.
\end{example}

\begin{example} 
\label{ex:7}
We end this section with an example of a non-compact Hamiltonian $S^1$-space with a {\bf non-proper} moment map {\bf with a critical point} for which the conclusion on the fundamental groups in Theorem~\ref{thm:2}  fails to hold. Let us consider $M=S^2\setminus{N} \times S^2\setminus{N}$ with symplectic form $2 \rho \oplus \rho$, where $\rho$ is a symplectic form on $S^2$ with total area $1$ and $N$ is the north pole, equipped with the standard diagonal circle action (cf. Figure~\ref{fig:3}). This action is Hamiltonian and its moment map $\phi$  has a unique critical value at $0$ corresponding to the fixed point $(S,S)$, where $S$ is the south pole of the sphere (see footnote~\ref{fn:1} in page~\pageref{fn:1}, for the moment map expression and compare Figures~\ref{fig:1} and \ref{fig:3}). We can see that the reduced spaces have different fundamental groups. Indeed, for $a\in (0,1)$, the reduced spaces $M_a=\phi^{-1}(a)/S^1$ are spheres while, for $a\in (2,3)$, they are spheres minus two points. That is,  $\pi_1(M_a)=\{1\}$ for $a\in(0,1)$, and  $\pi_1(M_a)= \bbZ$  for $a\in(2,3)$. 
\begin{figure}[h!]
\begin{center}
\epsfxsize=.4\textwidth
    \psfrag{H}{$\scriptsize{\phi(M)}$}
    \psfrag{0}{$\scriptsize{0}$} 
    \psfrag{1}{$\scriptsize{1}$}
    \psfrag{2}{$\scriptsize{2}$}
    \psfrag{3}{$\scriptsize{3}$}
    \psfrag{SS}{$\scriptsize{(S,S)}$}
    \leavevmode
    \epsfbox{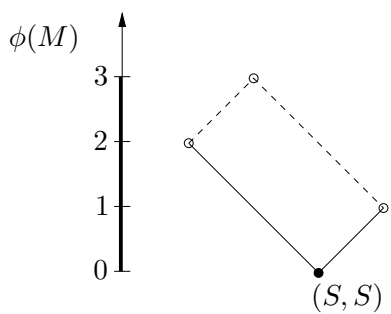}
\end{center}
\caption{Moment map image for an extended $T^2$-action on $S^2\setminus\{N\} \times S^2\setminus\{N\}$.}\label{fig:3}
\end{figure}
\end{example}


\begin{thebibliography}{SFO}

\bibitem[BHPV]{BHPV} W. Barth, K. Hulek, C. Peters and A. Van de Ven, {\em Compact Complex Surfaces}, Results in Mathematics and Related Areas. 3rd Series. A Series of Modern Surveys in Mathematics, {\bf 4}, Springer-Verlag, Berlin, (2004). 

\bibitem[BP]{BP} M. Brion and C. Procesi, {\em Action d'un tore dans une variet\'{e} projective}, Operator Algebras, Unitary Representations, Enveloping Algebras, and Invariant Theory, Progress in Mathematics, {\bf 93}, Birkh\"{a}user, Boston, (1991).

\bibitem[CGKZ]{CGKZ} D.J. Collins, R.I. Grigorchuk, P.F. Kurchanov and H. Zieschang, {\em Combinatorial Group Theory and Applications to Geometry}, Encyclopaedia of mathematical sciences, {\bf 58}, Springer, (1998).

\bibitem[FGM]{FGM} M. Fernandez, A. Gray and J.W. Morgan, {\em Compact symplectic manifolds with free circle actions and Massey products}, Michigan Math. J. {\bf 38} (1991), 271-283.

\bibitem[G]{G} L. Godinho, {\em Blowing up symplectic orbifolds}, Ann. Global Anal. Geom. {\bf 20} (2001), 117-162. 

\bibitem[L]{Le} E. Lerman, {\em Symplectic cuts}, Math. Res. Lett. {\bf 2} (1995), 247-258.

\bibitem[LMTW]{LMTW} E. Lerman, E. Meinrenken, S. Tolman and C. Woodward,  {\em Non abelian convexity by symplectic cuts}, Topology {\bf 37} (1998), 245-259.

\bibitem[Li]{L} H. Li, {\em $\pi_1$ of Hamiltonian $S^1$-actions}, Proc. Amer. Math. Soc. {\bf 131} (2003), 3579-3582.

\bibitem[K]{K} D. Kotschick, {\em Free circle actions with contractible orbits on symplectic manifolds}, SG$\backslash$0410243  (2004).

\bibitem[M]{M} T. Matumoto, {\em On diffeomorphisms of a K3 surface}, Algebraic and Topological Theories (Kinosaki, 1984), 616-621, Kinokuniya, Tokyo (1986).

\bibitem[MD1]{MD1} D. McDuff, {\em The moment map for circle actions on symplectic manifolds}, Journal of Geometry and Physics {\bf 5} (1998), no 2, 149-160.

\bibitem[MD2]{MD2} D. McDuff, {\em Examples of simply-connected symplectic non-K\"{a}hlerian manifolds}, Differential Geometry {\bf 20} (1984), 267-277.

\bibitem[Mi]{Mi} J. Milnor, {\em Morse Theory}, Princeton University Press, Princeton, N.J., (1963).

\bibitem[O]{O} K. Ono, {\em Obstruction to circle group actions preserving symplectic structure}, Hokkaido Math. Journal {\bf 21} (1992), 99-102.
 
\bibitem[OR]{OR} J.P. Ortega and T. Ratiu, {\em Momentum Maps and Hamiltonian Reduction}, Progress in Mathematics, {\bf 222}, Birkh\"{a}user, Boston, (2003).

\bibitem[R]{R} A. Ranicki, {\em High Dimensional Knot theory, Algebraic Surgery in Codimension 2}, Springer Monographs in Mathematics, {\bf 26}, Springer-Verlag, Berlin, (1998).

\end{thebibliography}
\end{document}